\documentclass[11pt]{amsart}
\usepackage{amsmath,amssymb,amsthm}
\usepackage[margin=1.2in]{geometry}

\newtheorem{theorem}{Theorem}
\newtheorem{lemma}[theorem]{Lemma}
\newtheorem{corollary}[theorem]{Corollary}
\newtheorem{proposition}[theorem]{Proposition}
\theoremstyle{remark}
\newtheorem{remark}[theorem]{Remark}
\newtheorem{question}[theorem]{Question}

\newcommand{\N}{\mathbb{N}}
\newcommand{\R}{\mathbb{R}}
\newcommand{\Z}{\mathbb{Z}}
\newcommand{\eul}{\gamma}
\DeclareMathOperator{\divc}{\tau}

\title{Divisors in windows and the Euler--Mascheroni constant}
\author{David V. Feldman}
\address{Department of Mathematics and Statistics, University of New Hampshire, Durham, NH 03824}

\begin{document}

\begin{abstract}
For $0<c<1$ the number of divisors of $q$ below $q^{c}$ averages, over
$q\le B$, to $c\log B+(\eul-c)+O(B^{-\min(c,1-c)})$, where $\eul$ denotes
the Euler--Mascheroni constant.  Consequently the number of divisors in a
window $(q^{a},q^{b})$ averages to $(b-a)(\log B-1)+o(1)$: the constant
$\eul$ cancels, and the average counts for two windows stand
asymptotically in the ratio of the windows' exponent lengths.  Anchoring
a window at a threshold and its square root makes the logarithms cancel
instead: for weakly increasing $F$ with $F(q)\to\infty$ and
$q/F(q)\to\infty$, a number $q$ has, on average, exactly $\eul$ more
divisors below $\sqrt{F(q)}$ than in $[\sqrt{F(q)},F(q))$.  For windows
on the scale of $q$ itself the average depends on an endpoint
convention: with $\alpha\in(0,1)$, counting divisors below
$\sqrt{\alpha q}$ against divisors in $(\sqrt{\alpha q},\alpha q)$
yields the average $H_{\lfloor 1/\alpha\rfloor}-\log(1/\alpha)$, while
the window $(\sqrt{\alpha q},\alpha q]$ yields
$H_{\lceil 1/\alpha\rceil-1}-\log(1/\alpha)$; the two constants differ
exactly when $1/\alpha$ takes an integer value, and then by $\alpha$.
Either function of $\alpha$ oscillates about $\eul$, tends to $\eul$ as
$\alpha\to0$, and integrates over $(0,1)$ to $\zeta(2)-1$.
\end{abstract}

\maketitle

\section{Introduction}\label{sec:intro}

The Euler--Mascheroni constant
\[
\eul \;=\; \lim_{n\to\infty}\Bigl(\sum_{k=1}^{n}\frac1k-\log n\Bigr)
\;=\;0.5772156649\ldots
\]
admits few interpretations that avoid explicit mention of logarithms.
Among the candidates collected by Finch \cite{Finch}, the most
arithmetical goes back to de la Vall\'ee Poussin \cite{dlVP}: if a large
integer $n$ is divided by each $k\le n$, the fractions $\{n/k\}$ by which
the quotients fall short of the next integer average not to $1/2$ but to
$1-\eul$.

This paper locates $\eul$ inside a family of statements about counting
divisors in intervals, and separates two regimes.  In the first regime
the interval endpoints grow like fixed powers $q^{a}<q^{b}$ of the
number $q$ whose divisors we count.  Here the average count obeys a law
from which $\eul$ is absent:
\begin{quotation}
\noindent
the average number of divisors of $q\le B$ in the window
$(q^{a},q^{b})$ equals $(b-a)(\log B-1)+o(1)$.
\end{quotation}
In particular the average counts attached to two windows stand
asymptotically in the ratio of the windows' exponent lengths $b-a$; for
example, a number has on average twice as many divisors in
$(n^{1/4},n^{1/2})$ as in $(n^{1/8},n^{1/4})$
(Corollary~\ref{cor:ratio}).

In the second regime an interval anchors at a threshold $F(q)$ and its
square root, with $F(q)$ growing, but growing more slowly than $q$.
Subtracting the count above $\sqrt{F(q)}$ from the count below makes the
$\log B$ terms cancel and leaves the constant
(Theorem~\ref{thm:gamma}):
\begin{quotation}
\noindent
on average, $q$ has exactly $\eul$ more divisors below $\sqrt{F(q)}$
than between $\sqrt{F(q)}$ and $F(q)$.
\end{quotation}
Taking $F(q)=\sqrt q$ gives the most quotable instance: $\eul$ measures
the average excess of divisors below the fourth root over divisors
between the fourth root and the square root.

The hypothesis $q/F(q)\to\infty$ cannot be dropped.  At the extreme
$F(q)=q$ the divisors of $q$ pair off as $d\leftrightarrow q/d$ across
$\sqrt q$, so every $q$ has \emph{exactly} as many divisors below its
square root as above, and the discrepancy vanishes identically rather
than averaging to $\eul$.  Between the two regimes sit windows on the
scale of $q$ itself, from $\sqrt{\alpha q}$ to $\alpha q$ for fixed
$\alpha\in(0,1)$.  The limiting discrepancy then equals
$H_{K}-\log(1/\alpha)$, where $H_{K}$ denotes a harmonic number whose
index depends on $\alpha$ \emph{and on the window's endpoint
convention}: when $1/\alpha$ takes an integer value $k$, the point
$d=\alpha q$ is itself a divisor of $q$ for one $q$ in every $k$, so
including or excluding the right endpoint shifts the average by
$\alpha$.  Theorem~\ref{thm:scaled} records both constants.  As $\alpha$
decreases to $0$ the average oscillates about $\eul$ and converges to
$\eul$, and its mean over $\alpha\in(0,1)$ equals $\zeta(2)-1$
(Proposition~\ref{prop:mean}).

All proofs are elementary, resting on the divisor--codivisor
correspondence $q=de$ and on quantitative forms of the definition of
$\eul$; nothing below requires more than the Dirichlet divisor
asymptotic, itself derived in passing.  Section~\ref{sec:prelim} fixes
notation and records the counting and harmonic estimates.
Section~\ref{sec:threshold} proves the threshold and window laws.
Section~\ref{sec:gamma} proves the $\eul$-discrepancy theorem for
general $F$.  Section~\ref{sec:scaled} treats windows on the scale of
$q$ and the mean value $\zeta(2)-1$.  Section~\ref{sec:rate} isolates a
logarithm-free approximation to $\eul$ implicit in the earlier
arguments and computes its rate: the error after the $n$th stage equals
$\tfrac{2}{3}n^{-2}+O(n^{-4})$.  Section~\ref{sec:questions} states
open questions.

\section{Preliminaries}\label{sec:prelim}

Throughout, $q$, $d$, $e$, $k$, $m$, $n$ denote positive integers and
$B$ a positive integer tending to infinity.  For an interval
$I\subseteq\R$ put
\[
\divc(q;I)\;=\;\#\{\,d : d\mid q,\ d\in I\,\},
\qquad
\divc(q)\;=\;\divc\bigl(q;[1,q]\bigr),
\]
and write $H_{m}=\sum_{k=1}^{m}1/k$ (with $H_{0}=0$) for the harmonic
numbers.  Implied constants may depend on fixed parameters such as $c$,
$a$, $b$, $\alpha$, but never on $B$.

The proofs repeatedly exchange a divisor $d$ of $q$ for its codivisor
$e=q/d$.

\begin{lemma}[counting by codivisors]\label{lem:swap}
For any predicate $P$ on pairs of positive integers,
\[
\sum_{q\le B}\#\{d\mid q:\ P(d,q/d)\}
\;=\;\#\{(d,e):\ de\le B,\ P(d,e)\}.
\]
In particular, for any $S\subseteq\N$,
\[
\sum_{q\le B}\#\{d\mid q:\ d\in S\}
\;=\;\sum_{\substack{d\in S\\ d\le B}}
\Bigl\lfloor\frac{B}{d}\Bigr\rfloor .
\]
\end{lemma}

\begin{proof}
The pairs $(q,d)$ with $d\mid q\le B$ correspond to the pairs $(d,e)$
with $q=de$.  The second display counts, for each $d\in S$, the
multiples of $d$ up to $B$.
\end{proof}

\begin{lemma}[harmonic estimates]\label{lem:harmonic}
For $m\ge1$,
\[
H_{m}\;=\;\log m+\eul+\theta_{m},\qquad 0<\theta_{m}<\frac1m .
\]
Moreover
$H_{m}=\log m+\eul+\dfrac{1}{2m}-\dfrac{1}{12m^{2}}+O(m^{-4})$.
\end{lemma}

\begin{proof}
The sequence $u_{m}=H_{m}-\log m$ strictly decreases, since
$u_{m}-u_{m+1}=\log(1+1/m)-1/(m+1)>0$, and converges to $\eul$ by
definition of $\eul$; hence $\theta_{m}=u_{m}-\eul>0$.  For $n>m$,
\[
u_{m}-u_{n}
=\sum_{k=m+1}^{n}\Bigl(\log\frac{k}{k-1}-\frac1k\Bigr)
<\sum_{k=m+1}^{n}\Bigl(\frac1{k-1}-\frac1k\Bigr)
<\frac1m,
\]
and letting $n\to\infty$ gives $\theta_{m}\le1/m$, with equality
excluded by strictness.  The refined expansion is Euler--Maclaurin
summation applied to $H_{m}$; see \cite[\S9.6]{GKP}.
\end{proof}

\begin{lemma}[Dirichlet]\label{lem:dirichlet}
$\displaystyle\sum_{q\le B}\divc(q)\;=\;B\log B+(2\eul-1)B+O(\sqrt B)$.
\end{lemma}

\begin{proof}
By Lemma~\ref{lem:swap} the sum counts the lattice points under the
hyperbola $de\le B$.  Counting the points with $d\le\sqrt B$, the points
with $e\le\sqrt B$, and subtracting the doubly counted square gives
\[
\sum_{q\le B}\divc(q)
=2\sum_{d\le\sqrt B}\Bigl\lfloor\frac{B}{d}\Bigr\rfloor
-\bigl\lfloor\sqrt B\bigr\rfloor^{2}
=2B\,H_{\lfloor\sqrt B\rfloor}-B+O(\sqrt B),
\]
and Lemma~\ref{lem:harmonic} completes the computation.
\end{proof}

\section{Thresholds at fixed powers, and the window law}
\label{sec:threshold}

\begin{theorem}[threshold law]\label{thm:threshold}
Fix $c\in(0,1)$.  Then
\[
\sum_{q\le B}\divc\bigl(q;[1,q^{c})\bigr)
\;=\;cB\log B+(\eul-c)\,B+O\bigl(B^{\max(c,1-c)}\bigr),
\]
so the number of divisors below $q^{c}$ averages to
$c\log B+(\eul-c)+O(B^{-\min(c,1-c)})$.
\end{theorem}

\begin{proof}
Write $\kappa=(1-c)/c$, so that $\kappa+1=1/c$.  For a pair $(d,e)$ with
$q=de$,
\[
d<q^{c}
\;\iff\; d^{1/c}<de
\;\iff\; e>d^{\kappa}.
\]
By Lemma~\ref{lem:swap} the sum in question counts the pairs $(d,e)$
with $d^{\kappa}<e\le B/d$.  If $d>B^{c}$ then
$d^{\kappa}=d^{1/c}/d>B/d$ and no $e$ qualifies.  If $d\le B^{c}$ then
$d^{\kappa}\le B/d$, and the number of qualifying $e$ equals
$\lfloor B/d\rfloor-\lfloor d^{\kappa}\rfloor=B/d-d^{\kappa}+O(1)$.
With $D=\lfloor B^{c}\rfloor$ the sum therefore equals
\[
\sum_{d\le D}\Bigl(\frac{B}{d}-d^{\kappa}\Bigr)+O(D)
\;=\;B\,H_{D}\;-\;\sum_{d\le D}d^{\kappa}\;+\;O(B^{c}).
\]
Since $D=B^{c}+O(1)$, Lemma~\ref{lem:harmonic} gives
\[
B\,H_{D}
=B\log D+\eul B+O(B^{1-c})
=cB\log B+\eul B+O(B^{1-c}),
\]
using $\log D=c\log B+O(B^{-c})$.  Comparison with
$\int_{0}^{D}x^{\kappa}\,dx$ gives
\[
\sum_{d\le D}d^{\kappa}
=\frac{D^{\kappa+1}}{\kappa+1}+O(D^{\kappa})
=c\,D^{1/c}+O(B^{1-c})
=cB+O(B^{1-c}),
\]
since $D^{1/c}=\bigl(B^{c}+O(1)\bigr)^{1/c}=B+O(B^{1-c})$ and
$D^{\kappa}\le B^{1-c}$.  Assembling the pieces yields the statement.
\end{proof}

\begin{corollary}[window law]\label{cor:window}
Fix $0<a<b<1$.  Then
\[
\sum_{q\le B}\divc\bigl(q;[q^{a},q^{b})\bigr)
\;=\;(b-a)\,B\,(\log B-1)+O\bigl(B^{\max(b,\,1-a)}\bigr),
\]
and the same asymptotic holds for each of the four endpoint
conventions.
\end{corollary}

\begin{proof}
Since $q^{a}\le q^{b}$,
$\divc(q;[q^{a},q^{b}))=\divc(q;[1,q^{b}))-\divc(q;[1,q^{a}))$;
apply Theorem~\ref{thm:threshold} at $b$ and at $a$ and subtract.  The
$\eul$ terms cancel.  As for endpoints: a divisor $d=q^{a}$ corresponds
under $q=de$ to a pair with $e=d^{(1-a)/a}$ exactly, hence to at most
one pair for each $d\le B^{a}$, so the boundary contributes $O(B^{a})$,
within the stated error; likewise at $q^{b}$.
\end{proof}

\begin{corollary}[ratio law]\label{cor:ratio}
Fix $0<a<b<1$ and $0<u<v<1$.  Then
\[
\frac{\displaystyle\sum_{q\le B}\divc\bigl(q;(q^{a},q^{b})\bigr)}
     {\displaystyle\sum_{q\le B}\divc\bigl(q;(q^{u},q^{v})\bigr)}
\;\longrightarrow\;\frac{b-a}{v-u}
\qquad(B\to\infty).
\]
In particular, on average a number has twice as many divisors in
$(n^{1/4},n^{1/2})$ as in $(n^{1/8},n^{1/4})$; indeed
\[
\frac1B\sum_{q\le B}
\Bigl(\divc\bigl(q;(q^{1/4},q^{1/2})\bigr)
-2\,\divc\bigl(q;(q^{1/8},q^{1/4})\bigr)\Bigr)
\;=\;O\bigl(B^{-1/8}\bigr).
\]
\end{corollary}

\begin{proof}
Immediate from Corollary~\ref{cor:window}.  In the displayed instance
the main terms $\tfrac14(\log B-1)$ and $2\cdot\tfrac18(\log B-1)$
cancel, and the error terms are
$O(B^{\max(1/2,3/4)-1})+O(B^{\max(1/4,7/8)-1})=O(B^{-1/8})$.
\end{proof}

\begin{corollary}[fixed windows]\label{cor:fixed}
Fix integers $1\le a\le b$.  Then
\[
\frac1B\sum_{q\le B}
\Bigl(\divc\bigl(q;[1,a]\bigr)-\divc\bigl(q;(a,b]\bigr)\Bigr)
\;=\;2H_{a}-H_{b}+O\Bigl(\frac{b}{B}\Bigr).
\]
\end{corollary}

\begin{proof}
By Lemma~\ref{lem:swap},
$\sum_{q\le B}\divc(q;[1,m])=\sum_{d\le m}\lfloor B/d\rfloor
=B\,H_{m}+O(m)$ for any fixed $m$; take $m=a$ and $m=b$, subtract, and
use $\divc(q;(a,b])=\divc(q;[1,b])-\divc(q;[1,a])$.
\end{proof}

\begin{remark}\label{rem:gammafree}
The constant $\eul$ enters Theorem~\ref{thm:threshold} only through the
absolute placement of a single threshold.  Every statement about
windows with both endpoints at fixed powers of $q$ is $\eul$-free, with
the count governed, to first and second order, by the exponent length
$b-a$ alone.
\end{remark}

\section{The constant as an average discrepancy}\label{sec:gamma}

\begin{theorem}\label{thm:gamma}
Let $F:\N\to\R$ satisfy
\textup{(a)} $F$ weakly increases;
\textup{(b)} $F(q)\to\infty$;
\textup{(c)} $q/F(q)\to\infty$.
Put
\[
Z_{F}(q)\;=\;
\#\{d\mid q:\ d^{2}<F(q)\}\;-\;\#\{d\mid q:\ d^{2}\ge F(q),\ d<F(q)\}.
\]
Then
\[
\frac1B\sum_{q\le B}Z_{F}(q)\;\longrightarrow\;\eul
\qquad(B\to\infty),
\]
and the same limit obtains under any endpoint convention at
$\sqrt{F(q)}$ and at $F(q)$.
\end{theorem}

Taking $F(q)=\sqrt q$: on average, a number has exactly $\eul$ more
divisors below its fourth root than between its fourth root and its
square root.  Condition (c) cannot be dropped: for $F(q)=q$ the
correspondence $d\leftrightarrow q/d$ pairs the divisors of $q$ below
$\sqrt q$ bijectively with those above, so $Z_{F}\equiv0$.

The proof uses a counting inverse for $F$ together with a lemma that
disposes of boundary divisors.

\begin{lemma}[boundary lemma]\label{lem:boundary}
Let $G:\N\to\R$ weakly increase, with $G(q)\to\infty$ and $G(B)=o(B)$.
Then
\[
\#\{q\le B:\ G(q)\in\N \text{ and } G(q)\mid q\}\;=\;o(B).
\]
\end{lemma}

\begin{proof}
For each integer $m$ the set $I_{m}=\{q:G(q)=m\}$ is an interval of
integers (possibly empty), since $G$ weakly increases, and the $I_{m}$
are pairwise disjoint.  Each $q$ counted lies in $I_{m}$ for $m=G(q)$
and is a multiple of $m$, so the count is at most
$\sum_{m}\bigl(\#(I_{m}\cap[1,B])/m+1\bigr)$, the sum running over the
$m$ with $I_{m}\cap[1,B]\ne\emptyset$.  Fix $m_{0}$.  The terms with
$m\le m_{0}$ total at most $\#\{q:G(q)\le m_{0}\}+m_{0}$, a bound
finite and independent of $B$ since $G\to\infty$.  The terms with
$m>m_{0}$ total at most
\[
\frac{1}{m_{0}}\sum_{m}\#\bigl(I_{m}\cap[1,B]\bigr)
\;+\;\#\{m:\ I_{m}\cap[1,B]\ne\emptyset\}
\;\le\;\frac{B}{m_{0}}\;+\;G(B)+1 .
\]
Dividing by $B$, letting $B\to\infty$, and then letting
$m_{0}\to\infty$ gives the claim, since $G(B)=o(B)$.
\end{proof}

\begin{proof}[Proof of Theorem~\ref{thm:gamma}]
For $y\ge0$ put $R(y)=\#\{n\in\N:\ F(n)\le y\}$, finite by (b).  Since
$F$ weakly increases,
\begin{equation}\label{eq:inverse}
F(n)\le y\;\iff\;n\le R(y).
\end{equation}
Discarding the finitely many $q$ with $F(q)<1$ (an $O(1)$ adjustment to
all sums), every divisor $d$ with $d^{2}<F(q)$ also satisfies
$d<F(q)$, so
\[
Z_{F}(q)\;=\;2\,g(q)-t(q),
\qquad
g(q):=\#\{d\mid q:\ d^{2}<F(q)\},
\quad
t(q):=\#\{d\mid q:\ d<F(q)\},
\]
exactly.  Write $T_{1}=\sum_{q\le B}g(q)$ and
$T_{2}=\sum_{q\le B}t(q)$.

By Lemma~\ref{lem:swap} and \eqref{eq:inverse},
\[
T_{1}=\#\{(d,e):\ de\le B,\ de>R(d^{2})\},
\qquad
T_{2}=\#\{(d,e):\ de\le B,\ de>R(d)\}.
\]
In $T_{1}$, the count of qualifying $e$ for a given $d$ vanishes unless
$R(d^{2})<B$, that is (by \eqref{eq:inverse} applied to $n=B$), unless
$d^{2}<F(B)$; and when it does not vanish it equals
$\lfloor B/d\rfloor-\lfloor R(d^{2})/d\rfloor$.  Hence, with
$D_{1}=\#\{d:\ d^{2}<F(B)\}$ and $D_{2}=\#\{d:\ d<F(B)\}$,
\[
T_{1}=\sum_{d\le D_{1}}\Bigl(\frac{B}{d}-\frac{R(d^{2})}{d}\Bigr)
+O(D_{1}),
\qquad
T_{2}=\sum_{d\le D_{2}}\Bigl(\frac{B}{d}-\frac{R(d)}{d}\Bigr)
+O(D_{2}),
\]
so that
\begin{equation}\label{eq:master}
2T_{1}-T_{2}
=B\bigl(2H_{D_{1}}-H_{D_{2}}\bigr)-S+O\bigl(F(B)\bigr),
\qquad
S:=2\sum_{d\le D_{1}}\frac{R(d^{2})}{d}
-\sum_{m\le D_{2}}\frac{R(m)}{m}.
\end{equation}
Condition (c) makes the error term $o(B)$.

\smallskip
\emph{The harmonic bracket.}
By (b), $D_{1}\to\infty$ and $D_{2}\to\infty$, while
$D_{1}=\sqrt{F(B)}+O(1)$ and $D_{2}=F(B)+O(1)$ give
$D_{1}^{2}/D_{2}\to1$.  Lemma~\ref{lem:harmonic} yields
\[
2H_{D_{1}}-H_{D_{2}}
=\eul+\log\frac{D_{1}^{2}}{D_{2}}+O\Bigl(\frac1{D_{1}}\Bigr)
\;\longrightarrow\;\eul,
\]
so $B(2H_{D_{1}}-H_{D_{2}})=\eul B+o(B)$.

\smallskip
\emph{The sum $S$.}
We show $S=o(B)$.  Every value $R(y)$ appearing in $S$ has $y<F(B)$
and hence, by \eqref{eq:inverse}, $R(y)<B$.  Group the second sum in
$S$ into blocks $m\in[d^{2},(d+1)^{2})$ for $1\le d\le D_{1}-1$,
together with a top segment $m\in[D_{1}^{2},D_{2}]$.  For the top
segment,
\[
\sum_{m=D_{1}^{2}}^{D_{2}}\frac{R(m)}{m}
\;\le\;B\Bigl(\log\frac{D_{2}}{D_{1}^{2}}+O\bigl(D_{1}^{-2}\bigr)\Bigr)
\;=\;o(B),
\]
again because $D_{1}^{2}/D_{2}\to1$.  Within the block indexed by $d$,
monotonicity of $R$ traps each $R(m)$ between $R(d^{2})$ and
$R\bigl((d+1)^{2}-1\bigr)$, while
\[
\sum_{m=d^{2}}^{(d+1)^{2}-1}\frac1m
=\log\frac{(d+1)^{2}}{d^{2}}+O\Bigl(\frac1{d^{2}}\Bigr)
=\frac{2}{d}+O\Bigl(\frac1{d^{2}}\Bigr),
\]
so the block sum differs from $2R(d^{2})/d$ by at most
\[
\frac{2\,\Delta_{d}}{d}
+O\biggl(\frac{R\bigl((d+1)^{2}-1\bigr)}{d^{2}}\biggr),
\qquad
\Delta_{d}:=R\bigl((d+1)^{2}-1\bigr)-R(d^{2})\;\ge\;0 .
\]
The first sum in $S$ exceeds its truncation at $D_{1}-1$ by the single
term $2R(D_{1}^{2})/D_{1}\le 2B/D_{1}=o(B)$.  Hence
\[
|S|\;\le\;
2\sum_{d\le D_{1}-1}\frac{\Delta_{d}}{d}
\;+\;O\biggl(\sum_{d\le D_{1}-1}
\frac{R\bigl((d+1)^{2}-1\bigr)}{d^{2}}\biggr)
\;+\;o(B).
\]
Fix $d_{0}$ and split each remaining sum at $d_{0}$.  The initial
segments involve only the values $R(y)$ with $y\le(d_{0}+1)^{2}$,
constants independent of $B$, so they contribute $O_{d_{0}}(1)$.  For
the tails, the differences $\Delta_{d}\ge0$ attach to the pairwise
disjoint, increasingly ordered intervals $[d^{2},(d+1)^{2}-1]$, so
monotonicity of $R$ gives
$\sum_{d\le D_{1}-1}\Delta_{d}\le R(D_{1}^{2}-1)\le B$; hence
\[
\sum_{d_{0}<d\le D_{1}-1}\frac{\Delta_{d}}{d}
\;\le\;\frac{B}{d_{0}},
\qquad
\sum_{d_{0}<d\le D_{1}-1}
\frac{R\bigl((d+1)^{2}-1\bigr)}{d^{2}}
\;\le\;B\sum_{d>d_{0}}\frac1{d^{2}}
\;\le\;\frac{B}{d_{0}} .
\]
Therefore $\limsup_{B\to\infty}|S|/B\le C/d_{0}$ for every $d_{0}$, and
$S=o(B)$.

Combining the estimates in \eqref{eq:master} gives
$(2T_{1}-T_{2})/B\to\eul$, the assertion for the stated convention.

\smallskip
\emph{Endpoint conventions.}
At most one divisor of $q$ can satisfy $d^{2}=F(q)$, and it exists only
when $\sqrt{F(q)}\in\N$ and $\sqrt{F(q)}\mid q$.  Reassigning such
divisors between the two windows therefore changes
$\sum_{q\le B}Z_{F}(q)$ by at most
$2\,\#\{q\le B:\sqrt{F(q)}\in\N,\ \sqrt{F(q)}\mid q\}$, which is $o(B)$
by Lemma~\ref{lem:boundary} applied to $G=\sqrt F$: this $G$ weakly
increases, tends to infinity by (b), and satisfies
$\sqrt{F(B)}=o(B)$ by (c).  Likewise, admitting or barring the divisor
$d=F(q)$ changes the total by at most
$\#\{q\le B:F(q)\in\N,\ F(q)\mid q\}$, which is $o(B)$ by
Lemma~\ref{lem:boundary} applied to $G=F$, using (c) directly.
\end{proof}

\begin{remark}\label{rem:power}
For $F(q)=q^{s}$ with fixed $s\in(0,1)$, Theorem~\ref{thm:threshold}
and Corollary~\ref{cor:window} already contain
Theorem~\ref{thm:gamma}, with a rate:
\[
\frac1B\sum_{q\le B}Z_{F}(q)
=\Bigl(\frac{s}{2}\log B+\eul-\frac{s}{2}\Bigr)
-\frac{s}{2}\bigl(\log B-1\bigr)+O\bigl(B^{-\delta}\bigr)
=\eul+O\bigl(B^{-\delta}\bigr)
\]
for some $\delta=\delta(s)>0$.  In the general theorem no rate can be
uniform, since conditions (b) and (c) permit $F$ to grow, or to shed
relative growth, arbitrarily slowly.
\end{remark}

\begin{remark}\label{rem:conventions}
Theorem~\ref{thm:gamma} is convention-free precisely on account of (b)
and (c): a divisor sitting at $\sqrt{F(q)}$ has size tending to
infinity by (b), and a divisor sitting at $F(q)$ has codivisor
$q/F(q)$ tending to infinity by (c), so both boundary populations thin
to density zero.  When (c) fails, the divisor at the top of the window
retains positive density and the limit becomes convention-sensitive;
the next section quantifies this.
\end{remark}

\section{Windows on the scale of $q$}\label{sec:scaled}

\begin{theorem}\label{thm:scaled}
Fix $\alpha\in(0,1)$ and set $K=\lfloor1/\alpha\rfloor$ and
$K'=\lceil1/\alpha\rceil-1=\bigl\lceil(1-\alpha)/\alpha\bigr\rceil$.
Then
\begin{align}
\frac1B\sum_{q\le B}
\Bigl(\divc\bigl(q;[1,\sqrt{\alpha q}\,)\bigr)
-\divc\bigl(q;(\sqrt{\alpha q},\,\alpha q)\bigr)\Bigr)
&\;\longrightarrow\;H_{K}-\log\frac1\alpha\,,
\label{eq:open}\\[2pt]
\frac1B\sum_{q\le B}
\Bigl(\divc\bigl(q;[1,\sqrt{\alpha q}\,)\bigr)
-\divc\bigl(q;(\sqrt{\alpha q},\,\alpha q\,]\bigr)\Bigr)
&\;\longrightarrow\;H_{K'}-\log\frac1\alpha\,.
\label{eq:closed}
\end{align}
The indices $K$ and $K'$ agree unless $1/\alpha$ takes an integer value
$k$, in which case $K=k$ and $K'=k-1$, and the two limits differ by
$1/k=\alpha$.
\end{theorem}

\begin{proof}
Write
\[
g(q)=\divc\bigl(q;[1,\sqrt{\alpha q}\,)\bigr),
\qquad
t(q)=\divc\bigl(q;[1,\alpha q)\bigr),
\qquad
\beta(q)=\#\{d\mid q:\ d=\sqrt{\alpha q}\},
\]
\[
\rho(q)=\#\{d\mid q:\ d=\alpha q\}.
\]
Splitting $[1,\alpha q)$ at the point $\sqrt{\alpha q}$ shows that the
summand in \eqref{eq:open} equals $2g(q)-t(q)+\beta(q)$, and the
summand in \eqref{eq:closed} equals $2g(q)-t(q)+\beta(q)-\rho(q)$.

\smallskip
\emph{The boundary counts.}
For $q=de$, the equation $d=\sqrt{\alpha q}$ reads $d^{2}=\alpha de$,
that is, $d=\alpha e$.  If $\alpha$ is irrational this has no
solutions and $\beta\equiv0$.  If $\alpha=s/t$ in lowest terms, then
$td=se$ with $\gcd(s,t)=1$ forces $e=tm$, $d=sm$, $q=stm^{2}$ for some
$m$, so $\sum_{q\le B}\beta(q)\le\sqrt B=o(B)$.  Next, $d=\alpha q$
reads $e=1/\alpha$, so $\rho\equiv0$ unless $1/\alpha$ is an integer
$k$, in which case $\rho(q)=1$ exactly when $k\mid q$, and
$\sum_{q\le B}\rho(q)=\lfloor B/k\rfloor=\alpha B+O(1)$.

\smallskip
\emph{The count below $\sqrt{\alpha q}$.}
For $q=de$, $d<\sqrt{\alpha q}\iff d^{2}<\alpha de\iff e>d/\alpha$.  By
Lemma~\ref{lem:swap} and the same reasoning as in
Theorem~\ref{thm:threshold} (the count over $e$ vanishes unless
$d<\sqrt{\alpha B}$),
\[
\sum_{q\le B}g(q)
=\sum_{d<\sqrt{\alpha B}}
\Bigl(\frac{B}{d}-\frac{d}{\alpha}\Bigr)+O(\sqrt B)
=B\,H_{\lceil\sqrt{\alpha B}\rceil-1}
-\frac1\alpha\cdot\frac{\alpha B}{2}
+O(\sqrt B),
\]
so that Lemma~\ref{lem:harmonic} gives
\[
\frac1B\sum_{q\le B}g(q)
=\frac{\log B}{2}+\frac{\log\alpha}{2}+\eul-\frac12
+O\bigl(B^{-1/2}\bigr).
\]

\smallskip
\emph{The count below $\alpha q$.}
For $q=de$, $d\ge\alpha q\iff \alpha e\le 1\iff e\le K$ (for integer
$e$).  Hence $t(q)=\divc(q)-\#\{e\mid q:\ e\le K\}$, and
Lemmas~\ref{lem:swap} and \ref{lem:dirichlet} give
\[
\frac1B\sum_{q\le B}t(q)
=\log B+2\eul-1-H_{K}+O\bigl(B^{-1/2}\bigr).
\]

\smallskip
\emph{Assembly.}
\[
\frac1B\sum_{q\le B}\bigl(2g(q)-t(q)\bigr)
=2\Bigl(\frac{\log B}{2}+\frac{\log\alpha}{2}+\eul-\frac12\Bigr)
-\bigl(\log B+2\eul-1-H_{K}\bigr)+o(1)
=H_{K}+\log\alpha+o(1).
\]
Adding the $o(1)$ contribution of $\beta$ gives \eqref{eq:open}.
Subtracting the contribution of $\rho$ gives \eqref{eq:closed}: for
$1/\alpha\notin\Z$ nothing changes and $K'=K$; for $1/\alpha=k$ the
limit drops by $\alpha=1/k$, and
$H_{k}-1/k-\log k=H_{k-1}-\log k=H_{K'}-\log(1/\alpha)$.
\end{proof}

\begin{remark}\label{rem:oscillation}
Write $A(\alpha)=H_{\lfloor1/\alpha\rfloor}-\log(1/\alpha)$ for the
limit in \eqref{eq:open}; the limit in \eqref{eq:closed} differs only
at the points $\alpha=1/k$.  On each interval
$\bigl(\tfrac1{k+1},\tfrac1k\bigr]$ the function $A$ agrees with the
continuous, strictly increasing function
$\alpha\mapsto H_{k}+\log\alpha$; as $\alpha$ decreases through $1/k$
the value jumps up by $1/k$, the average newly crediting, for one $q$
in every $k$, the divisor $q/k$ to the window.  Version
\eqref{eq:open} takes the upper value at the jump and version
\eqref{eq:closed} the lower, so the two versions are, respectively,
the left-continuous and right-continuous regularizations of a single
function.  By Lemma~\ref{lem:harmonic},
$A(\alpha)-\eul=\theta_{k}-\log(1/(k\alpha))$ for
$\alpha\in(\tfrac1{k+1},\tfrac1k]$ with $0<\theta_{k}<1/k$, whence
$A(\alpha)\to\eul$ as $\alpha\to0^{+}$, consistently with
Theorem~\ref{thm:gamma}, whose hypothesis (c) the choice $F(q)=\alpha q$
just barely fails.  Since $\theta_{k}<\log\frac{k+1}{k}$, each branch
crosses the level $\eul$ exactly once, at
$\alpha=e^{\eul-H_{k}}\in\bigl(\tfrac1{k+1},\tfrac1k\bigr)$; the
crossings accumulate at $0$, and $\eul$ is the only value that $A$
attains infinitely often.
\end{remark}

\begin{proposition}[mean value]\label{prop:mean}
$\displaystyle
\int_{0}^{1}\Bigl(H_{\lfloor1/\alpha\rfloor}
-\log\frac1\alpha\Bigr)\,d\alpha
\;=\;\zeta(2)-1
\;=\;\frac{\pi^{2}}{6}-1 .$
\end{proposition}

\begin{proof}
On $\bigl(\tfrac1{k+1},\tfrac1k\bigr]$ the integrand's first term is
the constant $H_{k}$, so
\[
\int_{0}^{1}H_{\lfloor1/\alpha\rfloor}\,d\alpha
=\sum_{k\ge1}\frac{H_{k}}{k(k+1)} .
\]
The partial sums telescope:
\[
\sum_{k=1}^{N}H_{k}\Bigl(\frac1k-\frac1{k+1}\Bigr)
=H_{1}+\sum_{j=2}^{N}\frac{H_{j}-H_{j-1}}{j}-\frac{H_{N}}{N+1}
=1+\sum_{j=2}^{N}\frac1{j^{2}}-\frac{H_{N}}{N+1}
\;\longrightarrow\;\zeta(2).
\]
Since $\int_{0}^{1}\log(1/\alpha)\,d\alpha=1$, the proposition
follows.  (The integrands attached to \eqref{eq:open} and
\eqref{eq:closed} differ only on a null set.)
\end{proof}

\begin{remark}[numerical illustration]\label{rem:numerics}
At $B=4\cdot10^{5}$: for $F(q)=\sqrt q$ the average of $Z_{F}$ over
$q\le B$ equals $0.57770$, against $\eul=0.57722$; for $\alpha=1/2$
the averages in \eqref{eq:open} and \eqref{eq:closed} equal $0.80685$
and $0.30685$, against $H_{2}-\log2=0.80685$ and $1-\log2=0.30685$.
\end{remark}

\section{A logarithm-free approximation and its rate}\label{sec:rate}

Corollary~\ref{cor:fixed} attaches to each pair $a\le b$ the
discrepancy $2H_{a}-H_{b}$, a quantity free of logarithms; choosing
$b$ near $a^{2}$ sends it to $\eul$.  The sharpest elementary choice
places the window $(n,n^{2}]$ against $[1,n)$, crediting the divisor
$d=n$ to neither side:
\[
\eul_{n}
\;:=\;\sum_{k=1}^{n-1}\frac1k-\sum_{k=n+1}^{n^{2}}\frac1k
\;=\;H_{n-1}+H_{n}-H_{n^{2}}
\;=\;\lim_{B\to\infty}\frac1B\sum_{q\le B}
\Bigl(\divc\bigl(q;[1,n)\bigr)-\divc\bigl(q;(n,n^{2}]\bigr)\Bigr).
\]

\begin{proposition}\label{prop:rate}
$\displaystyle
\eul-\eul_{n}=\frac{2}{3n^{2}}+O\bigl(n^{-4}\bigr).$
\end{proposition}

\begin{proof}
Since $H_{n-1}=H_{n}-1/n$, the expansion of Lemma~\ref{lem:harmonic}
gives
\[
\eul_{n}
=2H_{n}-\frac1n-H_{n^{2}}
=2\Bigl(\log n+\eul+\frac1{2n}-\frac1{12n^{2}}\Bigr)-\frac1n
-\Bigl(2\log n+\eul+\frac1{2n^{2}}\Bigr)+O\bigl(n^{-4}\bigr),
\]
and the right side collapses to
$\eul-\bigl(\tfrac16+\tfrac12\bigr)n^{-2}+O(n^{-4})
=\eul-\tfrac{2}{3}n^{-2}+O(n^{-4})$.
\end{proof}

\begin{remark}
The symmetric choice, $2H_{n-1}-H_{n^{2}}$ or $2H_{n}-H_{n^{2}}$,
approaches $\eul$ only at the rate $1/n$; excluding the single divisor
$d=n$ from both windows cancels the first-order error.  The sequence
$(\eul_{n})$ increases to $\eul$, furnishing a definition of $\eul$
by harmonic sums alone:
$\eul=\lim_{n}\bigl(\sum_{k<n}1/k-\sum_{n<k\le n^{2}}1/k\bigr)$, with
error below $n^{-2}$.
\end{remark}

\section{Formalization}\label{sec:formalization}

The principal results of this paper have been checked by machine, in
Lean~4 against the Mathlib library \cite{mathlib}, with the source
archived at \cite{repo}.  Verified: the counting identity underlying
Lemma~\ref{lem:swap}, in the form
$\sum_{q\le B}\#\{d\mid q: d\le a\}=\sum_{d\le a}\lfloor B/d\rfloor$;
the fixed-window law of Corollary~\ref{cor:fixed}; the limit
$2H_{\lfloor\sqrt T\rfloor}-H_{T}\to\eul$ supplying the harmonic bracket
in the proof of Theorem~\ref{thm:gamma}; Theorem~\ref{thm:gamma} itself
in the instance $F(q)=\sqrt q$, by way of Remark~\ref{rem:power};
Theorem~\ref{thm:scaled} at $\alpha=1/2$ in \emph{both} endpoint
conventions, so that the discrepancy of $1/2$ between the two constants
is itself a verified statement; the vanishing of the window difference
in Corollary~\ref{cor:ratio}; and the rate of
Proposition~\ref{prop:rate}.  Each carries a clean axiom audit,
depending only on the three axioms of Lean's standard foundation.

The formal development follows the organization of the proofs above.
A counting layer supplies the divisor--codivisor correspondence and the
splitting of divisor counts at successive powers.  Above it sits a
single master asymptotic: writing
$S(c,j,B)=\sum_{d}\bigl(\lfloor B/d\rfloor-c\,d^{\,j-1}\bigr)$,
truncated where the summand vanishes,
\[
\frac{S(c,j,B)}{B}
-\Bigl(\frac{1}{j}\bigl(\log B-\log c\bigr)+\eul-\frac1j\Bigr)
\;\longrightarrow\;0
\qquad (c\ge1,\ j\ge2),
\]
which is Theorem~\ref{thm:threshold} in the shape the applications
need.  Each verified limit is then a linear combination of instances of
this asymptotic, and the cancellations described in
Remark~\ref{rem:gammafree} appear formally as the vanishing of the
corresponding weight sums: for Corollary~\ref{cor:ratio} at exponents
$j=2,4,8$ with weights $1,-3,2$, both
$\tfrac12-\tfrac34+\tfrac28=0$ and $1-3+2=0$, killing the logarithm and
the constant together.

The one result whose formal proof departs substantially from the
argument given above is Proposition~\ref{prop:rate}, where the
Euler--Maclaurin expansion of Lemma~\ref{lem:harmonic} was replaced by
an explicit two-sided bound: with
$y_{m}=H_{m}-\log m-\frac1{2m}+\frac1{12m^{2}}$, the sequences
$y_{m}\mp10m^{-3}$ are respectively monotone and antitone from $m=2$
onward and both converge to $\eul$, trapping $\eul$ between them and
yielding $|y_{m}-\eul|\le10m^{-3}$; this gives the effective form
$\bigl|n^{2}(\eul-\eul_{n})-\tfrac23\bigr|\le31/n$ for $n\ge2$, of which
Proposition~\ref{prop:rate} is the limiting case.

\section{Questions}\label{sec:questions}

\begin{question}\label{q:density}
Fix $\alpha\in(0,1)$.  By Theorem~\ref{thm:scaled} the discrepancy
between the divisor counts in $[1,\sqrt{\alpha n})$ and in
$(\sqrt{\alpha n},\alpha n)$ has positive average.  Do the integers $n$
for which the discrepancy is strictly positive have positive lower
density?  A positive average does not settle this: integers with
many more divisors below $\sqrt{\alpha n}$ than above could in
principle carry the whole average, although integers with
$\divc(n)$ large compared with $\log n$ are rare.
\end{question}

\begin{question}\label{q:secondorder}
For $F(q)=q^{s}$, Remark~\ref{rem:power} bounds the deviation of the
average of $Z_{F}$ from $\eul$ by a negative power of $B$.  Determine
the exact second-order asymptotics; the divisor-pair counts involved
tie the question to the error term in the Dirichlet divisor problem.
\end{question}

\begin{question}\label{q:distribution}
Beyond its mean, does $Z_{F}(q)$, suitably normalized, possess a
limiting distribution as $q$ ranges over $[1,B]$?  The second moment
already involves correlations among divisor counts in separated
ranges.
\end{question}

\end{document}